\documentclass[12pt]{article}
\usepackage{amsmath}
\usepackage{amssymb}
\usepackage{stmaryrd}
\usepackage{amsthm}

\usepackage{cancel}

\usepackage[colorlinks=true,citecolor=black,linkcolor=black,urlcolor=blue]{hyperref}

\usepackage[dvipsnames]{pstricks} 
\usepackage{multido}

\usepackage{graphicx}
\usepackage{xcolor}

\newcommand{\showgrid}{}

\newcommand{\gridon}{\renewcommand{\showgrid}{\psset{subgriddiv=1,griddots=10,gridlabels=6pt}\psgrid}}

\gridon 

\newgray{hellgrau}{.89}

\usepackage[applemac]{inputenc}

\usepackage{algorithmic}

\newif\ifenglish

\englishtrue

\newif\ifvariant
\variantfalse

\def\bit{\begin{itemize}}
\def\eit{\end{itemize}}
\def\beq{\begin{equation}}
\def\eeq{\end{equation}}

\englishtrue

\def\of#1{\left(#1\right)} 
\def\pas#1{\left(#1\right)} 
\def\brk#1{\left[#1\right]} 

\def\defeq{\stackrel{\text{\tiny def}}{=}}

\def\F{{\mathbb F}}

\def\Z{{\mathbb Z}}

\def\0{{\mathbf 0}} 
\def\1{{\mathbf 1}} 

\def\CT1{CT1}
\def\xxxCT1{CT1}

\newtheorem{thm}{\ifenglish Theorem\else Satz\fi}[section] 

\newtheorem{pro}[thm]{Proposition}

\newtheorem{lem}[thm]{Lemma}
\newtheorem{rem}[thm]{\ifenglish Remark\else Bemerkung\fi}

\newtheoremstyle{excstyle}
  {1em}
  {2pt}
  {\sffamily\footnotesize\slshape}
  {0pt}
  {\sffamily\footnotesize\bfseries}
  {:}
  { }
  {}
\theoremstyle{excstyle}

\def\lemref#1{Lemma~\ref{#1}}

\def\proref#1{Proposition~\ref{#1}}

\def\figref#1{\ifenglish Figure\else Abbildung\fi~\ref{#1}}

\ifenglish\else\fi

\def\figref#1{Fi\-gu\-re~\ref{#1}}
\def\defeq{:=}
\def\qpoch#1#2#3{\pas{#1;#2}_{#3}}

\newgray{mfgray85}{0.85}

\newgray{mfgray70}{0.70}

\newgray{mfgray60}{0.60}

\title{The generating function of lozenge tilings for a ``quarter'' of a hexagon, obtained with non--intersecting lattice paths}

\author{Markus Fulmek\thanks{
Research supported by the National Research Network ``Analytic
Combinatorics and Probabilistic Number Theory'', funded by the
Austrian Science Foundation. 
}\\
\small Fakult\"at f\"ur Mathematik \\
\small Oskar-Morgenstern-Platz 1, A-1090 Wien, Austria \\
\small \tt Markus.Fulmek@Univie.Ac.At
}

\date{2020} 

\def\secA{\section}
\def\secB{\subsection}

\def\EM#1{{\em #1\/}}
\begin{document}
\bibliographystyle{plain}


\long\def\insertfigure#1#2{%
\begin{figure}%
\begin{center}%
\input graphics/#1%
\end{center}%
\caption{#2}%
\label{fig:#1}%
\end{figure}%
}

\long\def\inserttwocommentedfigures#1#2#3#4{%
\begin{figure}%
\begin{center}%
\input graphics/#1%
\hfil%
\input graphics/#2%
\end{center}%
{\small #4}%
\caption{#3}%
\label{fig:#1-#2}%
\end{figure}%
}

\long\def\insertcommentedfigure#1#2#3{%
\begin{figure}%
\begin{center}%
\input graphics/#1%
\end{center}%
{{\small #3}}
\caption{#2}%
\label{fig:#1}%
\end{figure}%
}

\parindent0pt
\parskip1em

\def\gf#1{{\mathfrak{gf}}\!\of{#1}}
\def\gfone#1{{\overline{\mathfrak{gf}}}\!\of{#1}}
\def\gftwo#1{{\widehat{\mathfrak{gf_1}}}\!\of{#1}}
\def\gfthree#1{{\widehat{\mathfrak{gf_0}}}\!\of{#1}}

\def\gfhalfeven#1{{\mathbf{gf}_0}\!\of{#1}}
\def\gfhalfodd#1{{\mathbf{gf}_1}\!\of{#1}}
\def\laiweight#1{\pas{X q^{#1}+Y q^{-\pas{#1}}}}
\def\maindet#1{{\mathbf d}\of{#1}}
\def\mfmat#1{{\mathbf m}\of{#1}}
\def\mfmatprime#1{{\mathbf m^\prime}\of{#1}}
\def\mfmatpprime#1{{\mathbf m^{\prime\prime}}\of{#1}}
\def\mfmatppprime{{\mathbf m^{\prime\prime\prime}}}
\def\mfsimple#1{{\mathbf s}\of{#1}}
\def\mffac#1{{\mathbf f}\of{#1}}
\def\mfdet#1{{\mathbf d}\of{#1}}
\def\minor#1#2#3#4{\brk{\maindet{#1}}_{#2,#3}^{#4}}
\def\mfminor#1#2#3#4{\brk{\mfdet{#1}}_{#2,#3}^{#4}}
\def\laplaceminor#1#2{\brk{\maindet{#1}}_{\cancel 1,\cancel #2}}

\def\mfentry#1{{\mathbf e}\of{#1}}

\def\prodcoeff#1#2#3{\mu^{#1}_{#2,#3}}

\def\mfseq{{\mathop {\mathbf a}}}

\def\qbinom#1#2{\genfrac{[}{]}{0pt}{}{#1}{#2}_q}
\def\qexp#1#2{e\of{#1,#2}}
\def\inner#1#2{\left\langle#1,#2\right\rangle}
\def\ratfield{\F}
\def\pprime{{\prime\prime}}
\def\mainminor#1#2#3{\left[{#1}\right]_{i\geq #2,j\geq #3}}
\def\step#1{\subsubsection*{#1:}}
\def\qop{E_q}
\def\qdiffop{D_q}
\def\id{I}

\def\thedet#1#2{{\mathbf d}\of{#1;#2}}
\def\mat#1{{\mathbf #1}}
\def\submat#1#2#3{{#1}_{\cancel{#2}\vert\cancel{#3}}}
\def\diag#1{{\mathbf{dg}}\of{#1}}

\newgray{backgroundgray}{0.965}

\maketitle

\begin{abstract}
In a recent preprint, Lai and Rohatgi 
compute the generating functions
of lozenge tilings of ``quartered hexagons with dents'' by applying the method of ``graphical condensation''.
The purpose of this note is to exhibit how
(a generalization of) Theorems 2.1 and 2.2 in Lai and Rohatgi's preprint 
can be achieved by the
Lindström--Gessel--Viennot method 
of non--intersecting lattice paths and a certain determinant evaluation.
\end{abstract}

\secA{Introduction}
In a recent preprint, Lai and Rohatgi \cite{Lai-Rohatgi:2020:TGFOHHAQH} compute the generating functions
of lozenge tilings of ``quartered hexagons with dents''; their method of proof is induction based on
``graphical condensation'' (i.e., an application of a certain Pfaffian identity to the enumeration of matchings).

The purpose of this note is to exhibit how (a generalization of) Theorems 2.1 and 2.2 in \cite{Lai-Rohatgi:2020:TGFOHHAQH}
can be achieved by the Lindström--Gessel--Viennot method \cite{Lindstroem:1973:OTVROIM,Gessel-Viennot:1998:DPAPP} of
non--intersecting lattice paths and the evaluation of the corresponding determinant. (Theorems 2.3 and 2.4
in \cite{Lai-Rohatgi:2020:TGFOHHAQH} most probably can be achieved in the same way, but we want to be brief here.)

This note is organized as follows: In section~\ref{sec:exposition}, we shall describe the generating functions
considered by Lai and Rohatgi, and show how the computation of these generating functions boils down (by the
well--known bijective correspondence between lozenge tilings and non--intersecting lattice paths
and the well--known Lindström--Gessel--Viennot argument) to
the evaluation of a certain determinant. In section~\ref{sec:evaluation}, we shall present two proofs for
the product formula giving the evaluation of this determinant.

\secA{Lozenge tilings of a ``quarter hexagon''}
\label{sec:exposition}

The literature on tilings enumerations is abundant (see,
for instance, \cite{Ciucu-Eisenkoelbl-Krattenthaler-Zare:2001:EOLT}), so we shall be brief.
For the experienced reader it certainly suffices to have a look at the left pictures of Figures~\ref{fig:lai22a-lai22aLGV}
and \ref{fig:lai22d-lai22dLGV},
which should make clear the concept of a ``quarter hexagon'' with ``dents'' (i.e., missing triangles) at its
base line. All \EM{vertical} lozenges of the tilings we consider here are \EM{labelled} with integers: This labelling
is \EM{vertically constant} and \EM{horizontally increasing by $1$ from
left to right}; in \figref{fig:lai22a-lai22aLGV}, this labelling 
starts at $1$, 
and in \figref{fig:lai22d-lai22dLGV}, this labelling starts at $0$.

Let $T$ be some lozenge tiling whose vertical
lozenges are labelled ${v_1, v_2, \dots, v_n}$, then the weight of $T$ is defined as
$$
w\of T\defeq\prod_{i=1}^n \frac{q^{v_i} + q^{-v_i}}2.
$$
Lai and Rohatgi \cite[Theorems 2.1--2.4]{Lai-Rohatgi:2020:TGFOHHAQH} computed the generating functions
of lozenge tilings for such ``quarter hexagons with dents'' of \EM{odd} or \EM{even} heights, and with
labels starting at $0$ or $1$ 
(see Figures~\ref{fig:lai22a-lai22aLGV}
and \ref{fig:lai22d-lai22dLGV}), which resulted in four theorems: The case where the labelling starts at $0$
(see \figref{fig:lai22d-lai22dLGV})
requires a slight modification of the weight function, but we shall only consider the other case (labelling starts at
$1$, see \figref{fig:lai22a-lai22aLGV}) and give a generalization which contains the cases of odd and even heights,
thus giving an alternative proof for Theorems 2.1 and 2.2 in \cite{Lai-Rohatgi:2020:TGFOHHAQH}.

%
\begin{figure}%
\begin{center}%
\psset{unit=0.6cm}
\begin{pspicture}(3.8,-0.2)(13.2,7.9942)
\pspolygon[linecolor=white,fillstyle=solid,fillcolor=backgroundgray](3.8,-0.2)(13.2,-0.2)(13.2,7.9942)(3.8,7.9942)
\psset{fillstyle=solid,linecolor=gray,linewidth=0.02}
\pspolygon[fillcolor=Tan](4.0,0.0)(4.5,0.86603)(5.5,0.86603)(5.0,0.0)
\psline[linewidth=0.07,linecolor=white](4.25,0.43301)(5.25,0.43301)
\pspolygon[fillcolor=Tan](5.0,0.0)(5.5,0.86603)(6.5,0.86603)(6.0,0.0)
\psline[linewidth=0.07,linecolor=white](5.25,0.43301)(6.25,0.43301)
\pspolygon[fillcolor=black](6.0,0.0)(6.5,0.86603)(7.0,0.0)
\pspolygon[fillcolor=black](7.0,0.0)(7.5,0.86603)(8.0,0.0)
\pspolygon[fillcolor=Mahogany](8.0,0.0)(7.5,0.86603)(8.5,0.86603)(9.0,0.0)
\pspolygon[fillcolor=black](9.0,0.0)(9.5,0.86603)(10.0,0.0)
\pspolygon[fillcolor=Mahogany](10.0,0.0)(9.5,0.86603)(10.5,0.86603)(11.0,0.0)
\pspolygon[fillcolor=black](11.0,0.0)(11.5,0.86603)(12.0,0.0)
\pspolygon[fillcolor=black](12.0,0.0)(12.5,0.86603)(13.0,0.0)
\pspolygon[fillcolor=Tan](4.5,0.86603)(5.0,1.7321)(6.0,1.7321)(5.5,0.86603)
\psline[linewidth=0.07,linecolor=white](4.75,1.299)(5.75,1.299)
\pspolygon[fillcolor=Tan](5.5,0.86603)(6.0,1.7321)(7.0,1.7321)(6.5,0.86603)
\psline[linewidth=0.07,linecolor=white](5.75,1.299)(6.75,1.299)
\pspolygon[fillcolor=Apricot](6.5,0.86603)(7.0,1.7321)(7.5,0.86603)(7.0,0.0)
\psline[linewidth=0.07,linecolor=white](6.75,1.299)(7.25,0.43301)
\rput[b](7.0,0.86603){ {\scriptsize\black $6$}}
\pspolygon[fillcolor=Tan](7.5,0.86603)(8.0,1.7321)(9.0,1.7321)(8.5,0.86603)
\psline[linewidth=0.07,linecolor=white](7.75,1.299)(8.75,1.299)
\pspolygon[fillcolor=Apricot](8.5,0.86603)(9.0,1.7321)(9.5,0.86603)(9.0,0.0)
\psline[linewidth=0.07,linecolor=white](8.75,1.299)(9.25,0.43301)
\rput[b](9.0,0.86603){ {\scriptsize\black $10$}}
\pspolygon[fillcolor=Tan](9.5,0.86603)(10.0,1.7321)(11.0,1.7321)(10.5,0.86603)
\psline[linewidth=0.07,linecolor=white](9.75,1.299)(10.75,1.299)
\pspolygon[fillcolor=Apricot](10.5,0.86603)(11.0,1.7321)(11.5,0.86603)(11.0,0.0)
\psline[linewidth=0.07,linecolor=white](10.75,1.299)(11.25,0.43301)
\rput[b](11.0,0.86603){ {\scriptsize\black $14$}}
\pspolygon[fillcolor=Apricot](11.5,0.86603)(12.0,1.7321)(12.5,0.86603)(12.0,0.0)
\psline[linewidth=0.07,linecolor=white](11.75,1.299)(12.25,0.43301)
\rput[b](12.0,0.86603){ {\scriptsize\black $16$}}
\pspolygon[fillcolor=Apricot](4.0,1.7321)(4.5,2.5981)(5.0,1.7321)(4.5,0.86603)
\psline[linewidth=0.07,linecolor=white](4.25,2.1651)(4.75,1.299)
\rput[b](4.5,1.7321){ {\scriptsize\black $1$}}
\pspolygon[fillcolor=Mahogany](5.0,1.7321)(4.5,2.5981)(5.5,2.5981)(6.0,1.7321)
\pspolygon[fillcolor=Mahogany](6.0,1.7321)(5.5,2.5981)(6.5,2.5981)(7.0,1.7321)
\pspolygon[fillcolor=Apricot](7.0,1.7321)(7.5,2.5981)(8.0,1.7321)(7.5,0.86603)
\psline[linewidth=0.07,linecolor=white](7.25,2.1651)(7.75,1.299)
\rput[b](7.5,1.7321){ {\scriptsize\black $7$}}
\pspolygon[fillcolor=Mahogany](8.0,1.7321)(7.5,2.5981)(8.5,2.5981)(9.0,1.7321)
\pspolygon[fillcolor=Apricot](9.0,1.7321)(9.5,2.5981)(10.0,1.7321)(9.5,0.86603)
\psline[linewidth=0.07,linecolor=white](9.25,2.1651)(9.75,1.299)
\rput[b](9.5,1.7321){ {\scriptsize\black $11$}}
\pspolygon[fillcolor=Mahogany](10.0,1.7321)(9.5,2.5981)(10.5,2.5981)(11.0,1.7321)
\pspolygon[fillcolor=Apricot](11.0,1.7321)(11.5,2.5981)(12.0,1.7321)(11.5,0.86603)
\psline[linewidth=0.07,linecolor=white](11.25,2.1651)(11.75,1.299)
\rput[b](11.5,1.7321){ {\scriptsize\black $15$}}
\pspolygon[fillcolor=Mahogany](4.5,2.5981)(4.0,3.4641)(5.0,3.4641)(5.5,2.5981)
\pspolygon[fillcolor=Tan](5.5,2.5981)(6.0,3.4641)(7.0,3.4641)(6.5,2.5981)
\psline[linewidth=0.07,linecolor=white](5.75,3.0311)(6.75,3.0311)
\pspolygon[fillcolor=Apricot](6.5,2.5981)(7.0,3.4641)(7.5,2.5981)(7.0,1.7321)
\psline[linewidth=0.07,linecolor=white](6.75,3.0311)(7.25,2.1651)
\rput[b](7.0,2.5981){ {\scriptsize\black $6$}}
\pspolygon[fillcolor=Tan](7.5,2.5981)(8.0,3.4641)(9.0,3.4641)(8.5,2.5981)
\psline[linewidth=0.07,linecolor=white](7.75,3.0311)(8.75,3.0311)
\pspolygon[fillcolor=Apricot](8.5,2.5981)(9.0,3.4641)(9.5,2.5981)(9.0,1.7321)
\psline[linewidth=0.07,linecolor=white](8.75,3.0311)(9.25,2.1651)
\rput[b](9.0,2.5981){ {\scriptsize\black $10$}}
\pspolygon[fillcolor=Mahogany](9.5,2.5981)(9.0,3.4641)(10.0,3.4641)(10.5,2.5981)
\pspolygon[fillcolor=Apricot](10.5,2.5981)(11.0,3.4641)(11.5,2.5981)(11.0,1.7321)
\psline[linewidth=0.07,linecolor=white](10.75,3.0311)(11.25,2.1651)
\rput[b](11.0,2.5981){ {\scriptsize\black $14$}}
\pspolygon[fillcolor=Tan](4.0,3.4641)(4.5,4.3301)(5.5,4.3301)(5.0,3.4641)
\psline[linewidth=0.07,linecolor=white](4.25,3.8971)(5.25,3.8971)
\pspolygon[fillcolor=Apricot](5.0,3.4641)(5.5,4.3301)(6.0,3.4641)(5.5,2.5981)
\psline[linewidth=0.07,linecolor=white](5.25,3.8971)(5.75,3.0311)
\rput[b](5.5,3.4641){ {\scriptsize\black $3$}}
\pspolygon[fillcolor=Mahogany](6.0,3.4641)(5.5,4.3301)(6.5,4.3301)(7.0,3.4641)
\pspolygon[fillcolor=Apricot](7.0,3.4641)(7.5,4.3301)(8.0,3.4641)(7.5,2.5981)
\psline[linewidth=0.07,linecolor=white](7.25,3.8971)(7.75,3.0311)
\rput[b](7.5,3.4641){ {\scriptsize\black $7$}}
\pspolygon[fillcolor=Mahogany](8.0,3.4641)(7.5,4.3301)(8.5,4.3301)(9.0,3.4641)
\pspolygon[fillcolor=Mahogany](9.0,3.4641)(8.5,4.3301)(9.5,4.3301)(10.0,3.4641)
\pspolygon[fillcolor=Apricot](10.0,3.4641)(10.5,4.3301)(11.0,3.4641)(10.5,2.5981)
\psline[linewidth=0.07,linecolor=white](10.25,3.8971)(10.75,3.0311)
\rput[b](10.5,3.4641){ {\scriptsize\black $13$}}
\pspolygon[fillcolor=Mahogany](4.5,4.3301)(4.0,5.1962)(5.0,5.1962)(5.5,4.3301)
\pspolygon[fillcolor=Mahogany](5.5,4.3301)(5.0,5.1962)(6.0,5.1962)(6.5,4.3301)
\pspolygon[fillcolor=Apricot](6.5,4.3301)(7.0,5.1962)(7.5,4.3301)(7.0,3.4641)
\psline[linewidth=0.07,linecolor=white](6.75,4.7631)(7.25,3.8971)
\rput[b](7.0,4.3301){ {\scriptsize\black $6$}}
\pspolygon[fillcolor=Mahogany](7.5,4.3301)(7.0,5.1962)(8.0,5.1962)(8.5,4.3301)
\pspolygon[fillcolor=Tan](8.5,4.3301)(9.0,5.1962)(10.0,5.1962)(9.5,4.3301)
\psline[linewidth=0.07,linecolor=white](8.75,4.7631)(9.75,4.7631)
\pspolygon[fillcolor=Apricot](9.5,4.3301)(10.0,5.1962)(10.5,4.3301)(10.0,3.4641)
\psline[linewidth=0.07,linecolor=white](9.75,4.7631)(10.25,3.8971)
\rput[b](10.0,4.3301){ {\scriptsize\black $12$}}
\pspolygon[fillcolor=Tan](4.0,5.1962)(4.5,6.0622)(5.5,6.0622)(5.0,5.1962)
\psline[linewidth=0.07,linecolor=white](4.25,5.6292)(5.25,5.6292)
\pspolygon[fillcolor=Tan](5.0,5.1962)(5.5,6.0622)(6.5,6.0622)(6.0,5.1962)
\psline[linewidth=0.07,linecolor=white](5.25,5.6292)(6.25,5.6292)
\pspolygon[fillcolor=Apricot](6.0,5.1962)(6.5,6.0622)(7.0,5.1962)(6.5,4.3301)
\psline[linewidth=0.07,linecolor=white](6.25,5.6292)(6.75,4.7631)
\rput[b](6.5,5.1962){ {\scriptsize\black $5$}}
\pspolygon[fillcolor=Tan](7.0,5.1962)(7.5,6.0622)(8.5,6.0622)(8.0,5.1962)
\psline[linewidth=0.07,linecolor=white](7.25,5.6292)(8.25,5.6292)
\pspolygon[fillcolor=Apricot](8.0,5.1962)(8.5,6.0622)(9.0,5.1962)(8.5,4.3301)
\psline[linewidth=0.07,linecolor=white](8.25,5.6292)(8.75,4.7631)
\rput[b](8.5,5.1962){ {\scriptsize\black $9$}}
\pspolygon[fillcolor=Mahogany](9.0,5.1962)(8.5,6.0622)(9.5,6.0622)(10.0,5.1962)
\pspolygon[fillcolor=Mahogany](4.5,6.0622)(4.0,6.9282)(5.0,6.9282)(5.5,6.0622)
\pspolygon[fillcolor=Mahogany](5.5,6.0622)(5.0,6.9282)(6.0,6.9282)(6.5,6.0622)
\pspolygon[fillcolor=Apricot](6.5,6.0622)(7.0,6.9282)(7.5,6.0622)(7.0,5.1962)
\psline[linewidth=0.07,linecolor=white](6.75,6.4952)(7.25,5.6292)
\rput[b](7.0,6.0622){ {\scriptsize\black $6$}}
\pspolygon[fillcolor=Mahogany](7.5,6.0622)(7.0,6.9282)(8.0,6.9282)(8.5,6.0622)
\pspolygon[fillcolor=Mahogany](8.5,6.0622)(8.0,6.9282)(9.0,6.9282)(9.5,6.0622)
\pspolygon[fillcolor=Tan](4.0,6.9282)(4.5,7.7942)(5.5,7.7942)(5.0,6.9282)
\psline[linewidth=0.07,linecolor=white](4.25,7.3612)(5.25,7.3612)
\pspolygon[fillcolor=Tan](5.0,6.9282)(5.5,7.7942)(6.5,7.7942)(6.0,6.9282)
\psline[linewidth=0.07,linecolor=white](5.25,7.3612)(6.25,7.3612)
\pspolygon[fillcolor=Apricot](6.0,6.9282)(6.5,7.7942)(7.0,6.9282)(6.5,6.0622)
\psline[linewidth=0.07,linecolor=white](6.25,7.3612)(6.75,6.4952)
\rput[b](6.5,6.9282){ {\scriptsize\black $5$}}
\pspolygon[fillcolor=Mahogany](7.0,6.9282)(6.5,7.7942)(7.5,7.7942)(8.0,6.9282)
\pspolygon[fillcolor=Mahogany](8.0,6.9282)(7.5,7.7942)(8.5,7.7942)(9.0,6.9282)
\end{pspicture}%
\hfil%
\psset{unit=0.6cm}
\begin{pspicture}(6.8,-0.7)(18.2,8.7)
\pspolygon[linecolor=white,fillstyle=solid,fillcolor=backgroundgray](6.8,-0.7)(18.2,-0.7)(18.2,8.7)(6.8,8.7)
\psset{fillstyle=none}
\psline[linestyle=dashed,linecolor=red](8.0,4.0)(17.0,8.5)
\psline[linecolor=lightgray](7.7,0.0)(17.5,0.0)
\rput(7.5,-0.2){ {\tiny 7}}
\rput(8.5,-0.2){ {\tiny 8}}
\rput(9.5,-0.2){ {\tiny 9}}
\rput(10.5,-0.2){ {\tiny 10}}
\rput(11.5,-0.2){ {\tiny 11}}
\rput(12.5,-0.2){ {\tiny 12}}
\rput(13.5,-0.2){ {\tiny 13}}
\rput(14.5,-0.2){ {\tiny 14}}
\rput(15.5,-0.2){ {\tiny 15}}
\rput(16.5,-0.2){ {\tiny 16}}
\psline[linecolor=lightgray](7.7,1.0)(17.5,1.0)
\rput(7.5,0.8){ {\tiny 5}}
\rput(8.5,0.8){ {\tiny 6}}
\rput(9.5,0.8){ {\tiny 7}}
\rput(10.5,0.8){ {\tiny 8}}
\rput(11.5,0.8){ {\tiny 9}}
\rput(12.5,0.8){ {\tiny 10}}
\rput(13.5,0.8){ {\tiny 11}}
\rput(14.5,0.8){ {\tiny 12}}
\rput(15.5,0.8){ {\tiny 13}}
\rput(16.5,0.8){ {\tiny 14}}
\psline[linecolor=lightgray](7.7,2.0)(17.5,2.0)
\rput(7.5,1.8){ {\tiny 3}}
\rput(8.5,1.8){ {\tiny 4}}
\rput(9.5,1.8){ {\tiny 5}}
\rput(10.5,1.8){ {\tiny 6}}
\rput(11.5,1.8){ {\tiny 7}}
\rput(12.5,1.8){ {\tiny 8}}
\rput(13.5,1.8){ {\tiny 9}}
\rput(14.5,1.8){ {\tiny 10}}
\rput(15.5,1.8){ {\tiny 11}}
\rput(16.5,1.8){ {\tiny 12}}
\psline[linecolor=lightgray](7.7,3.0)(17.5,3.0)
\rput(7.5,2.8){ {\tiny 1}}
\rput(8.5,2.8){ {\tiny 2}}
\rput(9.5,2.8){ {\tiny 3}}
\rput(10.5,2.8){ {\tiny 4}}
\rput(11.5,2.8){ {\tiny 5}}
\rput(12.5,2.8){ {\tiny 6}}
\rput(13.5,2.8){ {\tiny 7}}
\rput(14.5,2.8){ {\tiny 8}}
\rput(15.5,2.8){ {\tiny 9}}
\rput(16.5,2.8){ {\tiny 10}}
\psline[linecolor=lightgray](7.7,4.0)(17.5,4.0)
\rput(7.5,3.8){ {\tiny -1}}
\rput(8.5,3.8){ {\tiny 0}}
\rput(9.5,3.8){ {\tiny 1}}
\rput(10.5,3.8){ {\tiny 2}}
\rput(11.5,3.8){ {\tiny 3}}
\rput(12.5,3.8){ {\tiny 4}}
\rput(13.5,3.8){ {\tiny 5}}
\rput(14.5,3.8){ {\tiny 6}}
\rput(15.5,3.8){ {\tiny 7}}
\rput(16.5,3.8){ {\tiny 8}}
\psline[linecolor=lightgray](7.7,5.0)(17.5,5.0)
\rput(7.5,4.8){ {\tiny -3}}
\rput(8.5,4.8){ {\tiny -2}}
\rput(9.5,4.8){ {\tiny -1}}
\rput(10.5,4.8){ {\tiny 0}}
\rput(11.5,4.8){ {\tiny 1}}
\rput(12.5,4.8){ {\tiny 2}}
\rput(13.5,4.8){ {\tiny 3}}
\rput(14.5,4.8){ {\tiny 4}}
\rput(15.5,4.8){ {\tiny 5}}
\rput(16.5,4.8){ {\tiny 6}}
\psline[linecolor=lightgray](7.7,6.0)(17.5,6.0)
\rput(7.5,5.8){ {\tiny -5}}
\rput(8.5,5.8){ {\tiny -4}}
\rput(9.5,5.8){ {\tiny -3}}
\rput(10.5,5.8){ {\tiny -2}}
\rput(11.5,5.8){ {\tiny -1}}
\rput(12.5,5.8){ {\tiny 0}}
\rput(13.5,5.8){ {\tiny 1}}
\rput(14.5,5.8){ {\tiny 2}}
\rput(15.5,5.8){ {\tiny 3}}
\rput(16.5,5.8){ {\tiny 4}}
\psline[linecolor=lightgray](7.7,7.0)(17.5,7.0)
\rput(7.5,6.8){ {\tiny -7}}
\rput(8.5,6.8){ {\tiny -6}}
\rput(9.5,6.8){ {\tiny -5}}
\rput(10.5,6.8){ {\tiny -4}}
\rput(11.5,6.8){ {\tiny -3}}
\rput(12.5,6.8){ {\tiny -2}}
\rput(13.5,6.8){ {\tiny -1}}
\rput(14.5,6.8){ {\tiny 0}}
\rput(15.5,6.8){ {\tiny 1}}
\rput(16.5,6.8){ {\tiny 2}}
\psline[linecolor=lightgray](7.7,8.0)(17.5,8.0)
\rput(7.5,7.8){ {\tiny -9}}
\rput(8.5,7.8){ {\tiny -8}}
\rput(9.5,7.8){ {\tiny -7}}
\rput(10.5,7.8){ {\tiny -6}}
\rput(11.5,7.8){ {\tiny -5}}
\rput(12.5,7.8){ {\tiny -4}}
\rput(13.5,7.8){ {\tiny -3}}
\rput(14.5,7.8){ {\tiny -2}}
\rput(15.5,7.8){ {\tiny -1}}
\rput(16.5,7.8){ {\tiny 0}}
\psline[linecolor=lightgray](8.0,-0.3)(8.0,8.3)
\psline[linecolor=lightgray](9.0,-0.3)(9.0,8.3)
\psline[linecolor=lightgray](10.0,-0.3)(10.0,8.3)
\psline[linecolor=lightgray](11.0,-0.3)(11.0,8.3)
\psline[linecolor=lightgray](12.0,-0.3)(12.0,8.3)
\psline[linecolor=lightgray](13.0,-0.3)(13.0,8.3)
\psline[linecolor=lightgray](14.0,-0.3)(14.0,8.3)
\psline[linecolor=lightgray](15.0,-0.3)(15.0,8.3)
\psline[linecolor=lightgray](16.0,-0.3)(16.0,8.3)
\psline[linecolor=lightgray](17.0,-0.3)(17.0,8.3)
\psline[linewidth=0.07,linecolor=gray]{->}(7.0,4.0)(18.0,4.0)
\psline[linewidth=0.07,linecolor=gray]{->}(8.0,-0.5)(8.0,8.5)
\psline[linewidth=0.07,linecolor=blue,linearc=0.1](17.0,8.0)(17.0,7.0)(17.0,6.0)
\pscircle[linewidth=0.08,linecolor=blue,fillstyle=solid,fillcolor=white](17.0,8.0){0.18}
\pscircle[linewidth=0.08,linecolor=blue,fillstyle=solid,fillcolor=white](17.0,6.0){0.18}
\psline[linewidth=0.07,linecolor=blue,linearc=0.1](15.0,7.0)(16.0,7.0)(16.0,6.0)(16.0,5.0)(17.0,5.0)
\pscircle[linewidth=0.08,linecolor=blue,fillstyle=solid,fillcolor=white](15.0,7.0){0.18}
\pscircle[linewidth=0.08,linecolor=blue,fillstyle=solid,fillcolor=white](17.0,5.0){0.18}
\psline[linewidth=0.07,linecolor=blue,linearc=0.1](13.0,6.0)(13.0,5.0)(14.0,5.0)(14.0,4.0)(15.0,4.0)(16.0,4.0)(16.0,3.0)(17.0,3.0)
\pscircle[linewidth=0.08,linecolor=blue,fillstyle=solid,fillcolor=white](13.0,6.0){0.18}
\pscircle[linewidth=0.08,linecolor=blue,fillstyle=solid,fillcolor=white](17.0,3.0){0.18}
\psline[linewidth=0.07,linecolor=blue,linearc=0.1](11.0,5.0)(11.0,4.0)(11.0,3.0)(12.0,3.0)(13.0,3.0)(14.0,3.0)(14.0,2.0)(15.0,2.0)(16.0,2.0)(16.0,1.0)(17.0,1.0)
\pscircle[linewidth=0.08,linecolor=blue,fillstyle=solid,fillcolor=white](11.0,5.0){0.18}
\pscircle[linewidth=0.08,linecolor=blue,fillstyle=solid,fillcolor=white](17.0,1.0){0.18}
\psline[linewidth=0.07,linecolor=blue,linearc=0.1](9.0,4.0)(9.0,3.0)(9.0,2.0)(10.0,2.0)(11.0,2.0)(11.0,1.0)(12.0,1.0)(12.0,0.0)(13.0,0.0)(14.0,0.0)(15.0,0.0)(16.0,0.0)(17.0,0.0)
\pscircle[linewidth=0.08,linecolor=blue,fillstyle=solid,fillcolor=white](9.0,4.0){0.18}
\pscircle[linewidth=0.08,linecolor=blue,fillstyle=solid,fillcolor=white](17.0,0.0){0.18}
\uput[135.0](16.75,8.0){\red  $2y=x$}
\end{pspicture}%
\end{center}%
{\small %
The left picture shows a ``quarter hexagon'' of \EM{odd} height $9$ in the triangular lattice: The base line
has $5$ ``dents'' (i.e., missing triangles; indicated in the picture by black colour). The picture also shows a
\EM{lozenge tiling} of this ``quarter hexagon with dents'', where the three possible orientations of lozenges
({\bf\color{Mahogany} left--tilted}, {\bf\color{Apricot} right--tilted} and {\bf\color{Tan} vertical}) are indicated
by three different colours: All the \EM{vertical lozenges} are labelled with integers (this labelling is
\EM{constant} vertically and \EM{increasing by $1$} horizontally from left to right), and the \EM{weight} of this
particular tiling is
{\footnotesize $$
w_{1}\cdot w_3\cdot w_5^2\cdot w_6^4 \cdot w_7^2\cdot w_9\cdot w_{10}^2\cdot w_{11}\cdot w_{12}\cdot w_{13}\cdot w_{14}^2\cdot w_{15}\cdot w_{16},
$$}%
where $w_i\defeq\frac{q^i + q^{-i}}2$. %
The \EM{non--intersecting lattice paths} corresponding to this tiling are indicated by white lines in the
left picture; and the right picture shows a ``reflected, rotated and tilted'' version of these paths in the lattice $\Z\times\Z$,
where horizontal edges $\pas{a,b}\to\pas{a+1,b}$ are labelled $a-2b$. Clearly, this correspondence is a \EM{bijection}
(introduced here ``graphically'') between lozenge tilings and non--intersecting lattice paths, and this bijection is
\EM{weight--preserving} if we define the
weight of some family $P$ of non--intersecting lattice paths as the product of $w_i$, where $i$ runs over the
labels of all horizontal edges belonging to paths in $P$. %
}%
\caption{Pictures corresponding to Figure 2.2.a in Lai and Rohatgi's preprint.}%
\label{fig:lai22a-lai22aLGV}%
\end{figure}%

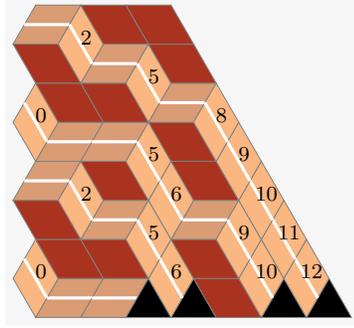
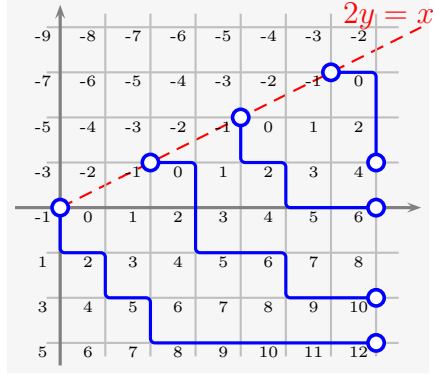
\begin{figure}%
\begin{center}%
\psset{unit=0.6cm}
\begin{pspicture}(3.3,-0.2)(11.2,7.1282)
\pspolygon[linecolor=white,fillstyle=solid,fillcolor=backgroundgray](3.3,-0.2)(11.2,-0.2)(11.2,7.1282)(3.3,7.1282)
\psset{fillstyle=solid,linecolor=gray,linewidth=0.02}
\pspolygon[fillcolor=Tan](4.0,0.0)(4.5,0.86603)(5.5,0.86603)(5.0,0.0)
\psline[linewidth=0.07,linecolor=white](4.25,0.43301)(5.25,0.43301)
\pspolygon[fillcolor=Tan](5.0,0.0)(5.5,0.86603)(6.5,0.86603)(6.0,0.0)
\psline[linewidth=0.07,linecolor=white](5.25,0.43301)(6.25,0.43301)
\pspolygon[fillcolor=black](6.0,0.0)(6.5,0.86603)(7.0,0.0)
\pspolygon[fillcolor=black](7.0,0.0)(7.5,0.86603)(8.0,0.0)
\pspolygon[fillcolor=Mahogany](8.0,0.0)(7.5,0.86603)(8.5,0.86603)(9.0,0.0)
\pspolygon[fillcolor=black](9.0,0.0)(9.5,0.86603)(10.0,0.0)
\pspolygon[fillcolor=black](10.0,0.0)(10.5,0.86603)(11.0,0.0)
\pspolygon[fillcolor=Apricot](3.5,0.86603)(4.0,1.7321)(4.5,0.86603)(4.0,0.0)
\psline[linewidth=0.07,linecolor=white](3.75,1.299)(4.25,0.43301)
\rput[b](4.0,0.86603){ {\scriptsize\black $0$}}
\pspolygon[fillcolor=Mahogany](4.5,0.86603)(4.0,1.7321)(5.0,1.7321)(5.5,0.86603)
\pspolygon[fillcolor=Mahogany](5.5,0.86603)(5.0,1.7321)(6.0,1.7321)(6.5,0.86603)
\pspolygon[fillcolor=Apricot](6.5,0.86603)(7.0,1.7321)(7.5,0.86603)(7.0,0.0)
\psline[linewidth=0.07,linecolor=white](6.75,1.299)(7.25,0.43301)
\rput[b](7.0,0.86603){ {\scriptsize\black $6$}}
\pspolygon[fillcolor=Mahogany](7.5,0.86603)(7.0,1.7321)(8.0,1.7321)(8.5,0.86603)
\pspolygon[fillcolor=Apricot](8.5,0.86603)(9.0,1.7321)(9.5,0.86603)(9.0,0.0)
\psline[linewidth=0.07,linecolor=white](8.75,1.299)(9.25,0.43301)
\rput[b](9.0,0.86603){ {\scriptsize\black $10$}}
\pspolygon[fillcolor=Apricot](9.5,0.86603)(10.0,1.7321)(10.5,0.86603)(10.0,0.0)
\psline[linewidth=0.07,linecolor=white](9.75,1.299)(10.25,0.43301)
\rput[b](10.0,0.86603){ {\scriptsize\black $12$}}
\pspolygon[fillcolor=Mahogany](4.0,1.7321)(3.5,2.5981)(4.5,2.5981)(5.0,1.7321)
\pspolygon[fillcolor=Tan](5.0,1.7321)(5.5,2.5981)(6.5,2.5981)(6.0,1.7321)
\psline[linewidth=0.07,linecolor=white](5.25,2.1651)(6.25,2.1651)
\pspolygon[fillcolor=Apricot](6.0,1.7321)(6.5,2.5981)(7.0,1.7321)(6.5,0.86603)
\psline[linewidth=0.07,linecolor=white](6.25,2.1651)(6.75,1.299)
\rput[b](6.5,1.7321){ {\scriptsize\black $5$}}
\pspolygon[fillcolor=Tan](7.0,1.7321)(7.5,2.5981)(8.5,2.5981)(8.0,1.7321)
\psline[linewidth=0.07,linecolor=white](7.25,2.1651)(8.25,2.1651)
\pspolygon[fillcolor=Apricot](8.0,1.7321)(8.5,2.5981)(9.0,1.7321)(8.5,0.86603)
\psline[linewidth=0.07,linecolor=white](8.25,2.1651)(8.75,1.299)
\rput[b](8.5,1.7321){ {\scriptsize\black $9$}}
\pspolygon[fillcolor=Apricot](9.0,1.7321)(9.5,2.5981)(10.0,1.7321)(9.5,0.86603)
\psline[linewidth=0.07,linecolor=white](9.25,2.1651)(9.75,1.299)
\rput[b](9.5,1.7321){ {\scriptsize\black $11$}}
\pspolygon[fillcolor=Tan](3.5,2.5981)(4.0,3.4641)(5.0,3.4641)(4.5,2.5981)
\psline[linewidth=0.07,linecolor=white](3.75,3.0311)(4.75,3.0311)
\pspolygon[fillcolor=Apricot](4.5,2.5981)(5.0,3.4641)(5.5,2.5981)(5.0,1.7321)
\psline[linewidth=0.07,linecolor=white](4.75,3.0311)(5.25,2.1651)
\rput[b](5.0,2.5981){ {\scriptsize\black $2$}}
\pspolygon[fillcolor=Mahogany](5.5,2.5981)(5.0,3.4641)(6.0,3.4641)(6.5,2.5981)
\pspolygon[fillcolor=Apricot](6.5,2.5981)(7.0,3.4641)(7.5,2.5981)(7.0,1.7321)
\psline[linewidth=0.07,linecolor=white](6.75,3.0311)(7.25,2.1651)
\rput[b](7.0,2.5981){ {\scriptsize\black $6$}}
\pspolygon[fillcolor=Mahogany](7.5,2.5981)(7.0,3.4641)(8.0,3.4641)(8.5,2.5981)
\pspolygon[fillcolor=Apricot](8.5,2.5981)(9.0,3.4641)(9.5,2.5981)(9.0,1.7321)
\psline[linewidth=0.07,linecolor=white](8.75,3.0311)(9.25,2.1651)
\rput[b](9.0,2.5981){ {\scriptsize\black $10$}}
\pspolygon[fillcolor=Tan](4.0,3.4641)(4.5,4.3301)(5.5,4.3301)(5.0,3.4641)
\psline[linewidth=0.07,linecolor=white](4.25,3.8971)(5.25,3.8971)
\pspolygon[fillcolor=Tan](5.0,3.4641)(5.5,4.3301)(6.5,4.3301)(6.0,3.4641)
\psline[linewidth=0.07,linecolor=white](5.25,3.8971)(6.25,3.8971)
\pspolygon[fillcolor=Apricot](6.0,3.4641)(6.5,4.3301)(7.0,3.4641)(6.5,2.5981)
\psline[linewidth=0.07,linecolor=white](6.25,3.8971)(6.75,3.0311)
\rput[b](6.5,3.4641){ {\scriptsize\black $5$}}
\pspolygon[fillcolor=Mahogany](7.0,3.4641)(6.5,4.3301)(7.5,4.3301)(8.0,3.4641)
\pspolygon[fillcolor=Apricot](8.0,3.4641)(8.5,4.3301)(9.0,3.4641)(8.5,2.5981)
\psline[linewidth=0.07,linecolor=white](8.25,3.8971)(8.75,3.0311)
\rput[b](8.5,3.4641){ {\scriptsize\black $9$}}
\pspolygon[fillcolor=Apricot](3.5,4.3301)(4.0,5.1962)(4.5,4.3301)(4.0,3.4641)
\psline[linewidth=0.07,linecolor=white](3.75,4.7631)(4.25,3.8971)
\rput[b](4.0,4.3301){ {\scriptsize\black $0$}}
\pspolygon[fillcolor=Mahogany](4.5,4.3301)(4.0,5.1962)(5.0,5.1962)(5.5,4.3301)
\pspolygon[fillcolor=Mahogany](5.5,4.3301)(5.0,5.1962)(6.0,5.1962)(6.5,4.3301)
\pspolygon[fillcolor=Tan](6.5,4.3301)(7.0,5.1962)(8.0,5.1962)(7.5,4.3301)
\psline[linewidth=0.07,linecolor=white](6.75,4.7631)(7.75,4.7631)
\pspolygon[fillcolor=Apricot](7.5,4.3301)(8.0,5.1962)(8.5,4.3301)(8.0,3.4641)
\psline[linewidth=0.07,linecolor=white](7.75,4.7631)(8.25,3.8971)
\rput[b](8.0,4.3301){ {\scriptsize\black $8$}}
\pspolygon[fillcolor=Mahogany](4.0,5.1962)(3.5,6.0622)(4.5,6.0622)(5.0,5.1962)
\pspolygon[fillcolor=Tan](5.0,5.1962)(5.5,6.0622)(6.5,6.0622)(6.0,5.1962)
\psline[linewidth=0.07,linecolor=white](5.25,5.6292)(6.25,5.6292)
\pspolygon[fillcolor=Apricot](6.0,5.1962)(6.5,6.0622)(7.0,5.1962)(6.5,4.3301)
\psline[linewidth=0.07,linecolor=white](6.25,5.6292)(6.75,4.7631)
\rput[b](6.5,5.1962){ {\scriptsize\black $5$}}
\pspolygon[fillcolor=Mahogany](7.0,5.1962)(6.5,6.0622)(7.5,6.0622)(8.0,5.1962)
\pspolygon[fillcolor=Tan](3.5,6.0622)(4.0,6.9282)(5.0,6.9282)(4.5,6.0622)
\psline[linewidth=0.07,linecolor=white](3.75,6.4952)(4.75,6.4952)
\pspolygon[fillcolor=Apricot](4.5,6.0622)(5.0,6.9282)(5.5,6.0622)(5.0,5.1962)
\psline[linewidth=0.07,linecolor=white](4.75,6.4952)(5.25,5.6292)
\rput[b](5.0,6.0622){ {\scriptsize\black $2$}}
\pspolygon[fillcolor=Mahogany](5.5,6.0622)(5.0,6.9282)(6.0,6.9282)(6.5,6.0622)
\pspolygon[fillcolor=Mahogany](6.5,6.0622)(6.0,6.9282)(7.0,6.9282)(7.5,6.0622)
\end{pspicture}%
\hfil%
\psset{unit=0.6cm}
\begin{pspicture}(4.8,-0.7)(14.2,7.7)
\pspolygon[linecolor=white,fillstyle=solid,fillcolor=backgroundgray](4.8,-0.7)(14.2,-0.7)(14.2,7.7)(4.8,7.7)
\psset{fillstyle=none}
\psline[linestyle=dashed,linecolor=red](6.0,3.0)(14.0,7.0)
\psline[linecolor=lightgray](5.7,0.0)(13.5,0.0)
\rput(5.5,-0.2){ {\tiny 5}}
\rput(6.5,-0.2){ {\tiny 6}}
\rput(7.5,-0.2){ {\tiny 7}}
\rput(8.5,-0.2){ {\tiny 8}}
\rput(9.5,-0.2){ {\tiny 9}}
\rput(10.5,-0.2){ {\tiny 10}}
\rput(11.5,-0.2){ {\tiny 11}}
\rput(12.5,-0.2){ {\tiny 12}}
\psline[linecolor=lightgray](5.7,1.0)(13.5,1.0)
\rput(5.5,0.8){ {\tiny 3}}
\rput(6.5,0.8){ {\tiny 4}}
\rput(7.5,0.8){ {\tiny 5}}
\rput(8.5,0.8){ {\tiny 6}}
\rput(9.5,0.8){ {\tiny 7}}
\rput(10.5,0.8){ {\tiny 8}}
\rput(11.5,0.8){ {\tiny 9}}
\rput(12.5,0.8){ {\tiny 10}}
\psline[linecolor=lightgray](5.7,2.0)(13.5,2.0)
\rput(5.5,1.8){ {\tiny 1}}
\rput(6.5,1.8){ {\tiny 2}}
\rput(7.5,1.8){ {\tiny 3}}
\rput(8.5,1.8){ {\tiny 4}}
\rput(9.5,1.8){ {\tiny 5}}
\rput(10.5,1.8){ {\tiny 6}}
\rput(11.5,1.8){ {\tiny 7}}
\rput(12.5,1.8){ {\tiny 8}}
\psline[linecolor=lightgray](5.7,3.0)(13.5,3.0)
\rput(5.5,2.8){ {\tiny -1}}
\rput(6.5,2.8){ {\tiny 0}}
\rput(7.5,2.8){ {\tiny 1}}
\rput(8.5,2.8){ {\tiny 2}}
\rput(9.5,2.8){ {\tiny 3}}
\rput(10.5,2.8){ {\tiny 4}}
\rput(11.5,2.8){ {\tiny 5}}
\rput(12.5,2.8){ {\tiny 6}}
\psline[linecolor=lightgray](5.7,4.0)(13.5,4.0)
\rput(5.5,3.8){ {\tiny -3}}
\rput(6.5,3.8){ {\tiny -2}}
\rput(7.5,3.8){ {\tiny -1}}
\rput(8.5,3.8){ {\tiny 0}}
\rput(9.5,3.8){ {\tiny 1}}
\rput(10.5,3.8){ {\tiny 2}}
\rput(11.5,3.8){ {\tiny 3}}
\rput(12.5,3.8){ {\tiny 4}}
\psline[linecolor=lightgray](5.7,5.0)(13.5,5.0)
\rput(5.5,4.8){ {\tiny -5}}
\rput(6.5,4.8){ {\tiny -4}}
\rput(7.5,4.8){ {\tiny -3}}
\rput(8.5,4.8){ {\tiny -2}}
\rput(9.5,4.8){ {\tiny -1}}
\rput(10.5,4.8){ {\tiny 0}}
\rput(11.5,4.8){ {\tiny 1}}
\rput(12.5,4.8){ {\tiny 2}}
\psline[linecolor=lightgray](5.7,6.0)(13.5,6.0)
\rput(5.5,5.8){ {\tiny -7}}
\rput(6.5,5.8){ {\tiny -6}}
\rput(7.5,5.8){ {\tiny -5}}
\rput(8.5,5.8){ {\tiny -4}}
\rput(9.5,5.8){ {\tiny -3}}
\rput(10.5,5.8){ {\tiny -2}}
\rput(11.5,5.8){ {\tiny -1}}
\rput(12.5,5.8){ {\tiny 0}}
\psline[linecolor=lightgray](5.7,7.0)(13.5,7.0)
\rput(5.5,6.8){ {\tiny -9}}
\rput(6.5,6.8){ {\tiny -8}}
\rput(7.5,6.8){ {\tiny -7}}
\rput(8.5,6.8){ {\tiny -6}}
\rput(9.5,6.8){ {\tiny -5}}
\rput(10.5,6.8){ {\tiny -4}}
\rput(11.5,6.8){ {\tiny -3}}
\rput(12.5,6.8){ {\tiny -2}}
\psline[linecolor=lightgray](6.0,-0.3)(6.0,7.3)
\psline[linecolor=lightgray](7.0,-0.3)(7.0,7.3)
\psline[linecolor=lightgray](8.0,-0.3)(8.0,7.3)
\psline[linecolor=lightgray](9.0,-0.3)(9.0,7.3)
\psline[linecolor=lightgray](10.0,-0.3)(10.0,7.3)
\psline[linecolor=lightgray](11.0,-0.3)(11.0,7.3)
\psline[linecolor=lightgray](12.0,-0.3)(12.0,7.3)
\psline[linecolor=lightgray](13.0,-0.3)(13.0,7.3)
\psline[linewidth=0.07,linecolor=gray]{->}(6.0,-0.5)(6.0,7.5)
\psline[linewidth=0.07,linecolor=gray]{->}(5.0,3.0)(14.0,3.0)
\psline[linewidth=0.07,linecolor=blue,linearc=0.1](12.0,6.0)(13.0,6.0)(13.0,5.0)(13.0,4.0)
\pscircle[linewidth=0.08,linecolor=blue,fillstyle=solid,fillcolor=white](12.0,6.0){0.18}
\pscircle[linewidth=0.08,linecolor=blue,fillstyle=solid,fillcolor=white](13.0,4.0){0.18}
\psline[linewidth=0.07,linecolor=blue,linearc=0.1](10.0,5.0)(10.0,4.0)(11.0,4.0)(11.0,3.0)(12.0,3.0)(13.0,3.0)
\pscircle[linewidth=0.08,linecolor=blue,fillstyle=solid,fillcolor=white](10.0,5.0){0.18}
\pscircle[linewidth=0.08,linecolor=blue,fillstyle=solid,fillcolor=white](13.0,3.0){0.18}
\psline[linewidth=0.07,linecolor=blue,linearc=0.1](8.0,4.0)(9.0,4.0)(9.0,3.0)(9.0,2.0)(10.0,2.0)(11.0,2.0)(11.0,1.0)(12.0,1.0)(13.0,1.0)
\pscircle[linewidth=0.08,linecolor=blue,fillstyle=solid,fillcolor=white](8.0,4.0){0.18}
\pscircle[linewidth=0.08,linecolor=blue,fillstyle=solid,fillcolor=white](13.0,1.0){0.18}
\psline[linewidth=0.07,linecolor=blue,linearc=0.1](6.0,3.0)(6.0,2.0)(7.0,2.0)(7.0,1.0)(8.0,1.0)(8.0,0.0)(9.0,0.0)(10.0,0.0)(11.0,0.0)(12.0,0.0)(13.0,0.0)
\pscircle[linewidth=0.08,linecolor=blue,fillstyle=solid,fillcolor=white](6.0,3.0){0.18}
\pscircle[linewidth=0.08,linecolor=blue,fillstyle=solid,fillcolor=white](13.0,0.0){0.18}
\uput[135.0](14.5,6.75){\red  $2y=x$}
\end{pspicture}%
\end{center}%
{\small %
The left picture shows a ``quarter hexagon with dents'' of \EM{even} height $8$ 
in the triangular lattice, and the right picture shows the corresponding family of non--intersecting lattice paths
in the lattice $\Z\times\Z$ (as explained in \figref{fig:lai22a-lai22aLGV}).
}%
\caption{Pictures corresponding to Figure 2.2.d in Lai and Rohatgi's preprint.}%
\label{fig:lai22d-lai22dLGV}%
\end{figure}%

\secB{The bijective correspondence between lo\-zen\-ge tilings and non--intersecting lattice paths}
The literature on the connection between lozenge tilings and non--intersecting lattice paths is abundant
(see, for instance, \cite[Section 5]{Ciucu-Eisenkoelbl-Krattenthaler-Zare:2001:EOLT});
for the experienced reader it certainly suffices to have
a look at the pictures in Figures~\ref{fig:lai22a-lai22aLGV}
and \ref{fig:lai22d-lai22dLGV}: It is easy to
see that there is a weight--preserving bijection between lozenge tilings and families of non--intersecting lattice paths in the 
lattice $\Z \times \Z$ with steps to the right and downwards, where steps to the right
from $\pas{a,b}$ to $\pas{a+1,b}$ are labelled $a-2b$ and thus have weight 
$$
\frac{q^{a-2b} + q^{2a-b}}{2}
$$
(and all downward steps have weight $1$). As usual, the weight of a lattice path is the product of all
the weights of steps it consists of.

Clearly, the generating function $\gf{a,b,c,d}$ of all lattice paths from initial point $\pas{a, b}$
to terminal point $\pas{c, d}$ is zero for $a>c$ or $b<d$. For $a<c$ and $b>d$, 
we claim
\begin{equation}
\label{eq:gf-lattice-path-general}
\gf{a,b,c,d} = 
\frac{
2^{a-c}q^{\frac{\pas{a-c}\pas{a+c-4d-1}}2}\qpoch{q^{2\pas{b-d+1}}}{q^2}{c-a}\qpoch{-q^{2\pas{a-b-d}}}{q^2}{c-a}
}{
\qpoch{q^2}{q^2}{c-a}
}.
\end{equation}
Here, we used the standard $q$--Pochhammer notation $\qpoch aq0 \defeq 1$ and
\begin{align}
\qpoch aqn &\defeq\prod_{j=0}^{n-1}\pas{1-a\cdot q^j}
\label{eq:qpochplus}
\\
\qpoch aq{-n} &\defeq\frac1{\qpoch{a\cdot q^{-1}}{q^{-1}}n}
\label{eq:qpochminus}
\end{align}
for integer $n>0$.

Equation \eqref{eq:gf-lattice-path-general} follows immediately by showing that it fulfils the obvious
\EM{recursion} for the generating function of such weighted lattice paths
\begin{align*}
\gf{a,b,a,d} &\equiv 1, \\ 
\gf{a,b,c,b} &= \prod_{i=a-2b}^{c-2b-1}\frac{q^i + q ^{-i}}2
=2^{a-c}q^{\frac{\pas{a-c}\pas{a-4b+c-1}}2}\qpoch{-q^{2a-4b}}{q^2}{c-a}, \\ 
\gf{a,b,c,d} &= \frac{q^{a-2b} + q ^{2b-a}}2\,\gf{a+1,b,c,d} + \gf{a,b-1,c,d}
\end{align*}
for $a\leq c$ and $b\geq d$.
 
We have to specialize \eqref{eq:gf-lattice-path-general} to our situation, i.e., to initial points
$\pas{2i-1,i-1}$ and terminal points $\pas{2m-1+k,a_j}$
(see the right picture in \figref{fig:lai22a-lai22aLGV}), where $k=0$ or $k=1$ for quarter hexagons of \EM{odd} or \EM{even}
height, respectively, and $\pas{a_1,\dots,a_m}$ is a strictly increasing sequence of integers (where
the $a_j$ correspond to the positions of the ``dents'' in the lozenge tiling):
\begin{multline}
\label{eq:gf-lattice-path-special}
\gf{2i-1,i-1,2m-1+k,a_j} \\= 
\frac{
2^{2i-k-2m}q^{\frac{\pas{2i-k-2m}\pas{2i+k+2m-4a_j-3}}2}\qpoch{q^{4\pas{i-a_j}}}{q^4}{2m+k-2i}
}{
\qpoch{q^2}{q^2}{2m+k-2i}
}.
\end{multline}
(Note that this is zero for $a_j\geq i$.)

By the well--known the Lindström--Gessel--Viennot argument \cite{Lindstroem:1973:OTVROIM,Gessel-Viennot:1998:DPAPP},
the generating function of all families of non--intersecting lattice paths (which, by the weight--preserving
bijection sketched above is equal to the generating function of all lozenge tilings) can be written as the determinant
$$
\det_{1\leq i,j\leq m}\of{\gf{2i-1,i-1,2m-1+k,a_j}}.
$$
Theorems 2.1 and 2.2 in \cite{Lai-Rohatgi:2020:TGFOHHAQH} state that this determinant factorizes
completely for $k=0$ and $k=1$: Here, we shall show that this is true for general $k\geq 0$.

By the \EM{multilinearity} of the determinant, we may pull out the denominator in \eqref{eq:gf-lattice-path-special}
from all rows,
\bit
\item and powers of $q$ depending only on $i$ from the rows
\item and powers of $q$ depending only on $j$ from the columns,
\eit
which operation leaves the determinant
\begin{equation}
\label{eq:remaining-det}
\det_{1\leq i,j\leq m}\of{q^{-4 i a_j}\cdot\qpoch{q^{4\pas{i-a_j}}}{q^4}{2m+k-2i}}.
\end{equation}
Clearly, the claimed complete factorization of generating functions follows once we can
show that the above determinant factorizes: For simplicity's sake, we may
substitute $q^4\to q$ in \eqref{eq:remaining-det}.

\secA{The evaluation of the determinant}
\label{sec:evaluation}
It is not so easy to guess the correct evaluation of \eqref{eq:remaining-det} from the 
``known special cases'' Theorems 2.1 and 2.2 in \cite{Lai-Rohatgi:2020:TGFOHHAQH} and from computer
experiments, but once this is achieved, the proof is simple:
\begin{pro}
\label{pro:det1}
Let $\mfseq=\pas{a_1,a_2,\dots,a_m}$ be a strictly increasing sequence of $m$ integers, and let $k$ be a
nonnegative integer. Consider the $m\times m$--matrix
\begin{equation}
\label{eq:thematrix}
\mat m\of{k;\mfseq} = \pas{\frac{\qpoch{q^{i-a_j}}{q}{2m+k-2i}}{q^{i\, a_j}}}_{i,j=1}^m.
\end{equation}
Then we have
\begin{multline}
\label{eq:det1}
\det\pas{\mat m\of{k;\mfseq}}
=
q^{\frac{m\pas{m-1}\pas{2m+k-1}}2-\sum_{l=1}^m a_l\pas{2m-l}}\\
\times\pas{\prod_{j=1}^m\qpoch{q^{m-a_j}}qk}\prod_{i=1}^{m-1}\prod_{j=i+1}^m\pas{1-q^{a_i-a_j}}\pas{1-q^{-2m-k+1+a_i+a_j}}.
\end{multline}
\end{pro}

As a prerequisite for the proof of \proref{pro:det1}, recall
Dodgson's condensation formula \cite{Dodgson:1866}:
Assume $\mat m$ is an $\pas{m\times m}$--matrix, and write $\submat{\mat m}{i_1,\dots}{j_1,\dots}$
for the
submatrix obtained from $\mat m$ by deleting rows $\pas{i_1,\dots}$ and columns $\pas{j_1,\dots}$. Then for $m\geq 2$,
Dodgson's condensation formula states
\begin{multline}
\label{eq:dodgson}
\det\of{\mat m}\cdot\det\of{\submat{\mat m}{1,m}{1,m}} \\
=
\det\of{\submat{\mat m}{1}{1}}
\cdot
\det\of{\submat{\mat m}{m}{m}}
-
\det\of{\submat{\mat m}{1}{m}}
\cdot
\det\of{\submat{\mat m}{m}{1}}.
\end{multline}
(By convention, the determinant of a $\pas{0\times0}$--matrix equals $1$.)

\begin{rem}
Dodgson's condensation formula  is also known as
Desnanot--Jacobi's adjoint matrix theorem, see \cite[Theorem 3.12]{bressoud:1999}:
According to \cite{bressoud:1999}, Lagrange discovered this theorem for $n=3$,
Desnanot proved it for $n\leq 6$
and Jacobi published the general theorem \cite{jacobi:1841a},
see also \cite[vol.~I, pp.~142]{muir:1906}).
\end{rem}
\begin{proof}
Observe that for $\mat m = \mat m\of{k;\mfseq}$, all the submatrices appearing in Dodgson's condensation
formula \eqref{eq:dodgson} are of a
``similar type'', which we can describe in a simple manner by introducing the notation
$\diag{x_1,\dots,x_k}$ for the $k\times k$ diagonal matrix with entries $x_i$ on the main diagonal, and
defining
\begin{align*}
{}^\backprime\mfseq &= \pas{a_2,\dots,a_m}, \\
\mfseq^\prime &= \pas{a_1,\dots,a_{m-1}}, \\
{}^\backprime\mfseq^\prime &= \pas{a_2,\dots,a_{m-1}}, \\
\mfseq -1 &=\pas{a_1-1,\dots,a_m-1}.
\end{align*}
for any sequence $\mfseq = \pas{a_1,\dots,a_m}$ of length $m\geq 2$:
{\small
\begin{align*}
\submat{{\mat m}\of{k;\mfseq}}{1}{1} &= 
\diag{{q^{-1},\dots,q^{-\pas{m-1}}}}\cdot
\mat m\of{k;{}^\backprime\mfseq - 1}\cdot
\diag{{q^{-a_2},\dots,q^{-{a_m}}}}, \\
\submat{{\mat m}\of{k;\mfseq}}{n}{n} &= \mat m\of{k+2;\mfseq^\prime}, \\
\submat{{\mat m}\of{k;\mfseq}}{1}{n} &= 
\diag{{q^{-1},\dots,q^{-\pas{m-1}}}}\cdot
\mat m\of{k;\mfseq^\prime - 1}\cdot
\diag{{q^{-a_1},\dots,q^{-{a_{m-1}}}}}, \\
\submat{{\mat m}\of{k;\mfseq}}{n}{1} &= \mat m\of{k+2;{}^\backprime\mfseq}, \\
\submat{{\mat m}\of{k;\mfseq}}{1,n}{1,n} &= 
\diag{{q^{-2},\dots,q^{-\pas{m-1}}}}\cdot
\mat m\of{k+2;{}^\backprime\mfseq^\prime - 1}\cdot
\diag{{q^{-a_2},\dots,q^{-{a_{m-1}}}}}.
\end{align*}
}

Observe that both sides of \eqref{eq:det1} are zero iff $a_m\geq m$. So we may assume $a_m<m$:
But then the last element $a_{m-1}-1$ of ${}^\backprime\mfseq^\prime - 1$ is less than $m-2$, too,
and hence $\det\of{\submat{{\mat m}\of{k+2;\mfseq}}{1,n}{1,n}}\neq 0$.
Now writing $\thedet{k}{\mfseq}\defeq\det\pas{{\mat m}\of{k;\mfseq}}$,
by Dodgson's condensation formula \eqref{eq:dodgson} and the multiplicativity of the determinant
we see
{\small
\begin{multline}
\label{eq:dodgson-induction}
\thedet{k}{\mfseq}= \\ q^{-m+1}\cdot
\frac{
	q^{-a_m}\cdot\thedet{k+2}{\mfseq^\prime}\cdot\thedet{k}{{}^\backprime\mfseq - 1}-
	q^{-a_1}\cdot\thedet{k+2}{{}^\backprime\mfseq}\cdot\thedet{k}{\mfseq^\prime-1}
}{
	\thedet{k+2}{{}^\backprime\mfseq^\prime-1}
}.
\end{multline}
}

So we can use induction on the length $m$ of the sequence $\mfseq$: For $m=0$ and $m=1$, \eqref{eq:det1}
is obviously true; and for the inductive step, we may use \eqref{eq:dodgson-induction} (leading to a
lengthy, but straightforward computation).
\end{proof}

\secB{A second proof}
When I showed the determinant \eqref{eq:det1} 
to Christian Krattenthaler, he almost immediately
recognized it as a special case of his Lemma~4 in \cite{Krattenthaler:2001:ADC}, which we shall
repeat here for reader's convenience:
\begin{lem} 
\label{lem:Krat2}
Let $X_1,X_2,\dots,X_m,A_2,\dots,A_m$ be
indeterminates. Then there holds
{\small
\begin{multline} \label{eq:Krat2}
\det_{1\le i,j\le m} \big(( {C}/ {X_i}+A_m)( {C}/ {X_i}+A_{m-1})\dotsb( {C}/
{X_i}+A_{j+1})\\
\cdot (X_i+A_m)(X_i+A_{m-1})\dotsb(X_i+A_{j+1})\big)
=\prod _{i=2} ^{m}A_i^{i-1}\prod _{1\le i<j\le m}
(X_i-X_j)(1-C/X_iX_j).
\end{multline}
}
\end{lem}
Using Krattenthaler's Lemma 
has the big advantage that it \EM{yields} the determinant
evaluation (i.e., there is no need to \EM{guess} it).
\begin{proof}[Second proof of \proref{pro:det1}.]
By writing
$$
\qpoch{q^{i-a_j}}q{2m+k-2i} =
\qpoch{q^{i-a_j}}q{m-i}\cdot
\qpoch{q^{m-a_j}}q{k}\cdot
\qpoch{q^{m+k-a_j}}q{m-i}
$$
and pulling out the middle factors $\qpoch{q^{m-a_j}}q{k}$ from all columns of $\mat m\of{k;\mfseq}$, we obtain
\begin{align*}
\thedet{k}{\mfseq}
&=\pas{\prod_{j=1}^m\qpoch{q^{m-a_j}}qk}\det_{1\leq i,j\leq m}\pas{q^{-i a_j} \qpoch{q^{i-a_j}}q{m-i}\qpoch{q^{m+k-a_j}}q{m-i}}.
\end{align*}
By setting $X_j=q^{-a_j}$, expanding the $q$--Pochhammer symbols and pulling out powers of $q$, the determinant
from the above expression becomes
\begin{align*}
&\hphantom{=}\det_{1\leq i,j\leq m}\pas{X_j^i \qpoch{q^i X_j}q{m-i} \qpoch{q^{m+k} X_j}q{m-i}} \\
&= \det_{1\leq i,j\leq m}\Biggl(
	X_j^i
	\pas{\prod_{l=0}^{m-i-1}q^{i+l}\pas{q^{-i-l}-X_j}} \\
&\times
	\pas{\prod_{l=0}^{m-i-1}X_j q^{2m+k-1}\pas{X_j^{-1}q^{-2m-k+1}-q^{-m+l+1}}}
\Biggr) \\
&= \prod_{l=1}^m \pas{X_l^{m} q^{\binom{m}2 - \binom{l}2 + \pas{2m+k-1}\pas{m-l}}}\\
&\times\det_{1\leq i,j\leq m}\pas{
	\pas{\prod_{l=0}^{m-i-1}\pas{q^{-i-l}-X_j}}
	\pas{\prod_{l=0}^{m-i-1}\pas{X_j^{-1}q^{-2m-k+1}-q^{-m+l+1}}}
}.
\end{align*}

Now apply \lemref{lem:Krat2} with $X_j=q^{-aj}$, $A_i=q^{-i+1}$ and $C=q^{-2m-k+1}$
to obtain \eqref{eq:det1} (after some straightforward computation).
\end{proof}

\section*{Acknowledgement}
I am grateful to Christian Krattenthaler for helpful discussions.

\bibliography{/Users/mfulmek/Informations/TeX/database}

\end{document}